\documentclass{article}
\usepackage{amsmath,latexsym,amssymb,amsthm,enumerate,amscd}
\theoremstyle{plain}
\newtheorem{theorem}{Theorem}[section]

\newtheorem{corollary}[theorem]{Corollary}
\newtheorem{lemma}[theorem]{Lemma}
\theoremstyle{definition}

\newtheorem{remark}[theorem]{Remark}
\newtheorem{example}[theorem]{Example}

\newtheorem{question}[theorem]{Question}

\newtheorem{thm}{Theorem}

\newtheorem{rem}[thm]{Remark}
\newtheorem*{ucla}{Claim}

\newtheorem*{uthm}{Theorem}

\newtheorem*{unote}{Note}
\newtheorem*{ack}{Acknowledgment}

\numberwithin{equation}{section}
\numberwithin{table}{section} %TTT
\def\cha{\mathrm{char}\ }

\def\Hom{\mathrm{Hom}}

\def\Ann{\mathrm{Ann}\ }

\def\LEVALG{\mathrm{LevAlg}}
\def\PGOR{\mathrm{PGOR}}
\def\Pgl{\mathrm{PGL}}

\def\Soc{\mathrm{Soc}}

\providecommand{\bysame}{\makebox[3em]{\hrulefill}\thinspace}

\def\<{\left<}
\def\>{\right>}

\def\ns{\footnotesize \it}

\def\Z{\mathfrak{Z}}

\title{Hilbert functions of Gorenstein algebras associated to a
pencil of forms}
\author{Anthony Iarrobino\\[.05in]
{\ns Department of Mathematics, Northeastern University, Boston, MA 02115, USA.
}\\[.2in]}

\date{December 17, 2004}

\begin{document}

\maketitle
%\subjclass{Primary: 13D40; Secondary 14C05}
%\footnote\ddag{AMS 2000 Subj. Class, Primary: 14C05; Secondary: 13D40} 
\begin{abstract}Let $R$ be a polynomial ring in $r$ variables over an
infinite field
$K$, and denote by
$\mathcal{D}$ a corresponding dual ring, upon which $R$ acts as differential
operators. We study type two graded level Artinian algebras $A=R/I$,
having socle degree $j$. For each such algebra
$A$, we consider the family of Artinian Gorenstein [AG] quotients of $A$
having the same socle degree. \par 
By Macaulay duality, $A$ corresponds to a unique
2-dimensional vector space
$W_A$ of forms in
$\mathcal{D}_j$, and each such AG quotient of $A$ corresponds 
to a form in $W_A$ - up
to non-zero multiple. For
$W_A=\langle F,G\rangle$, each such AG quotient $A_\lambda$ corresponds
to an element of the pencil of forms (one dimensional subspaces) of $W_A$:
given
$F_\lambda=F+\lambda G, \lambda \in K\cup \infty$ we have  $A_\lambda =
R/{Ann}(F_\lambda)$. Our main result is a lower bound for the Hilbert function
$H(A_{\lambda_{gen}})$ of the generic Gorenstein quotient, in
terms of
$H(A)$, and the pair $H_F=H(R/\Ann F)$ and $H_G=H(R/\Ann G)$. This result
restricts the possible sequences $H$ that may occur as the Hilbert
function $H(A)$ for a type two level algebra $A$.
\end{abstract}

\section{Introduction}
Let  $R=k[x_1,\ldots ,x_r]$ be the polynomial ring in $r$ variables, over an infinite
field
$K$. We will assume also for simplicity of exposition that $\cha K=0$, but
all statements may be extended suitably to characteristic $p$  (see Remark
\ref{charp}).  We will consider only graded Artinian quotients
$A=R/I$ of $R$, and we denote by
$m$ the irrelevant maximal ideal
$m=(x_1,\ldots ,x_r)$. We denote by
$A_i$ the
$i$-th graded component of $A$. Recall that the \emph{socle} $\Soc (A)$ of
$A$ satisfies
\begin{equation}
\Soc(A)=(0:m)=\{f\in A\mid mf=0\},
\end{equation}
and the \emph{type} of $A$ is $t(A)=\dim_K \Soc(A)$. We will denote by
$j=j(A)$ the
\emph{socle degree} of
$A$, the largest integer such that $A_j\not= 0$, but $A_{j+1}=0$. An
Artinian algebra $A=R/I$ of socle degree $j$ is \emph{level}, if any of the
following equivalent conditions hold
\begin{enumerate}[i.]
\item
$\Soc(A)=A_j$,
\item The canonical module $\Hom(A,K)$ of $A$ is generated in a
single degree,
\item  Each $I_i, 0\le i\le j$ can be recovered from $I_j$ as follows:
\begin{equation}
\text { for } 0<i\le j \quad I_i=I_j:R_{j-i}=\{f\in
R_i\mid R_{j-i}\cdot f\subset I_j\}.
\end{equation}
\end{enumerate} 
Recall that $A$ is \emph{Artinian Gorenstein} if $A$ is level of type
one. Our main objects of study here will be type two level algebras $A$
and their Artinian Gorenstein [AG] quotients. 
\subsection{Recent results on level algebras}
We first briefly recall some recent work on graded level algebras. First,
the Artinian Gorenstein algebras have been the object of much study. For
height three, the structure theorem of Buchsbaum-Eisenbud
\cite{BE} has led  not only to a characterization of the Hilbert functions 
that may occur, but also to the irreducibility
\cite{Di} and smoothness \cite{Kle2} of the family
$\PGOR(H)$ parametrizing Artinian Gorenstein quotients of $R, r=3$, of a
given (symmetric) Hilbert function $H$. A second proof of smoothness follows from
results of M. Boij and A. Conca-G. Valla  (see \cite{Bj4,CoVa}, and 
\cite[\S 4.4]{IK} for a survey of related topics). One line of study relates
punctual subschemes $\Z$ of $\mathbb P^{r-1}$ to Gorenstein Artinian quotients of
their coordinate rings $\mathcal O_\Z$ \cite{Bj2,Ge,EmI,IK,Kle3}. Classical
apolarity, or the inverse systems of Macaulay provide
 a  connection to  sums of powers of linear forms, and to a classical
Waring problem for forms (see \cite{Ter,EmI,Ge,IK}). \par
When $r\ge 4$ $\PGOR(H)$ often has several irreducible components, a fact
first noted by M. Boij \cite{Bj3}, and elaborated by others (see
\cite{IS, Kle3}). The set of Gorenstein sequences --- ones that occur as Hilbert
functions of Artinian Gorenstein quotients of $R$ --- is not known for $r\ge 4$; for $r\ge
5$ they include non-unimodal sequences, with several maxima.
\par
 Level algebras $A$ of types $t(A)>1$ are a natural next
topic of study after the Gorenstein algebras, particularly in the low
embedding dimensions $r\le 3$ or even $r=4$, where the families of
Artinian Gorenstein quotients are better understood. When $r=2$, the family
$\LEVALG(H)$ is well understood (see \cite{I3,ChGe}): these
families are smooth, of known dimension, and their closures are the union
$\bigcup_{H'\le H} \LEVALG (H')$ of similar strata for termwise no greater
Hilbert functions of the same socle degree.  When
$r=3$, and $t=2$, tables of possible $H$ for small socle degree $j$,
possible resolutions, and many methods that are more general are given in
\cite{GHMS}; this case is also studied by F. Zanello \cite{Za1}. However the
possible sequences $H$ are not known when $r\ge 3$ even in the case $t=2$; and
although specialists believe there should be cases where $\LEVALG(H)$ for
$r=3,t=2$ have several irreducible components, this problem is still open. There
has been work connecting these results with the simultaneous Waring problem for
binary and ternary forms
\cite{Car,CaCh,I3}.  \par Several authors have studied the extremal Hilbert
functions for level algebras of given embedding dimension $r$, type $t$, and
socle degree $j$
\cite{BiGe,ChoI}. The minimal resolutions for compressed level algebras 
(those having maximum Hilbert function given $(r,t,j)$)  are
studied in
\cite{Bj1,MMN,Za1}. Also Zanello has obtained results about extremal Hilbert
functions for level algebras,  given the pair $H_{j-1}$ and $H_j=2$
\cite{Za2}. For $r\ge 3$ there is much to be learned about families of level
algebras of given Hilbert functions, even when
$t=2$.
Certainly, pencils of curves on $\mathbb P^2$ have been long a topic of geometric
study; however, a stronger connection (but see \cite{ChGe}) has yet to be made
between on the one hand the traditional geometrical approach to pencils of curves
and their singularities, and on the other hand the study of the level
algebras associated to these pencils.
\par
In this article we show some inequalities connecting the Hilbert function of type
two level algebras, and the Hilbert function of their Artinian Gorenstein
quotients. These results had been embedded in the longer
preprint \cite{I4}, which will now be refocussed on refinements of the
numerical results and on parametrization.\par In section \ref{subsinvsys} we give
notation, and briefly state the main results, and in section \ref{subsquest} we present
further context, including the questions that motivated us. In Section
\ref{sec2} we prove our results and give examples.
\subsection{Inverse systems}\label{subsinvsys}
Let $\mathcal{D}=K[X_1,\ldots ,X_r]$ denote a second polynomial ring. The ring $R$
acts on $\mathcal{D}$ as partial differential operators: for $ h\in R, F\in
\mathcal{D}$
\begin{equation}  h\circ
F=h(\partial/\partial X_1,\ldots, \partial/\partial X_r)\circ F.
\end{equation} 
The pairing 
\begin{equation}
\sigma_j:R_j\times \mathcal{D}_j\to K, \quad \sigma_j(h,F)=h\circ F,
\end{equation}
 is exact. This is the \emph{apolarity} or \emph{Macaulay
duality} action of $R$ on $\mathcal{D}$.  A type $t$ level algebra $A=R/I$
 of socle
degree $j$ corresponds via the Macaulay duality to a unique
$t$-dimensional vector space $W_A$,
\begin{equation}
W_A= \{ F\in \mathcal{D}_j\mid  I\circ F=0\}= \{ F\in \mathcal{D}_j\mid 
I_j\circ F=0\}.
\end{equation}
Thus $W_A$ is the perpendicular space to $I_j$ in the exact duality
between $R_j$ and $\mathcal{D}_j$, and $R\circ W_A$, may be regarded
as the dualizing module $\hat{A}=\Hom(A,K)$ to $A$. The Hilbert function $H(A)$
is the sequence $H(A)_i=\dim_K A_i$. We
have, for
 $0\le u\le
j$,
\begin{align}\label{eqduala}
R_u\circ W_A&=I_{j-u}^\perp\subset \mathcal D_{j-u}, \text { and }\\
H(A)_{j-u}&=\dim_K R_u\circ W_A.\label{eqdualb}
\end{align}
\begin{remark}\label{remgen} A one-dimensional subspace $E\subset W_A$
corresponds to an Artinian Gorenstein [AG] quotient $R/\Ann E$ of the level
algebra $A$, having the same socle degree $j$ as $A$. We parametrize these spaces $E$
as points of the projective space ${\mathbb P}(W_A)$ associated
to $W_A$. The Hilbert function
$H_E=H(R/\Ann E)$ is evidently semicontinuous: $H_E>T$ termwise for some fixed
sequence $T=(t_0,\ldots ,t_j)$ defines an open subset of ${\mathbb P}(W_A)$ since $\dim
R_u\circ E>t_{j-u}$ is an open condition. Thus, among the
Hilbert functions of Gorenstein quotients of $A$ having the same socle degree
$j$, there is a termwise maximum 
$H(E_{gen})$, that occurs for $E$ belonging to an open dense subset of
the projective variety
${\mathbb P}(W_A)$.
\par When the type of $A$ is two, then
$W_A=\langle F,G\rangle$ is two-dimensional; the one-dimensional subspaces
constitute a \emph{pencil} of forms $F_\lambda=F+\lambda G,
\lambda \in K\cup
\infty={\mathbb P}^1_K={\mathbb P}(W_A)$ (here we set $F_\infty =G$).
 Each AG quotient algebra
$A_\lambda =R/\Ann(F+\lambda G)$ has the same socle degree $j=j(A)$ as
$A$, and these comprise all the Gorenstein quotients of $A$
having socle degree $j(A)$. Thus, the family $A_\lambda, \lambda \in
{\mathbb P}^1$ is the pencil of Artinian Gorenstein quotients associated to the
pencil
$F+\lambda G$. We let $A_F=R/\Ann F$ and $A_G=R/\Ann G$, and set $H_F=H(A_F),
H_G=H(A_G)$.
\end{remark}
 We focus here on the type two case, and on the pencil of
degree-$j$ homogeneous forms or hypersurfaces,
$F_\lambda \in\mathcal{D}_j$ and their symmetric Hilbert functions
$H_\lambda=H(A_\lambda)$. By \eqref{eqdualb} $(H_\lambda)_{j-u}$ is the dimension of the
space of order $u$ partial derivatives of $F_\lambda$. Evidently, the set of Hilbert
functions
$H(A_\lambda)$ that occur is a
$\Pgl(r-1)$ invariant of the level algebra $A$. Also, Remark \ref{remgen}
implies that there is a termwise maximum value $H(A_{\lambda_{gen}})$, that
occurs for an open dense set of
$\lambda\in
\mathbb P^1$. 
\par
We now state the most important part of our main result, Theorem \ref{mainthm}.
For
$i=j-u$, we let
\begin{equation*}
d_i(F,G)= \dim_K \langle R_u\circ F\rangle\cap \langle R_u\circ
G\rangle = (H_F)_i+(H_G)_i-H(A)_i,
\end{equation*}
 be the overlap dimension, satisfying $d_i(F,G)=H(R/(\Ann F+\Ann G))_i$ (see equation
\eqref{eqdi}ff).
\begin{uthm}\label{umainthm}
 Let $A=R/\Ann (F,G)$ be a type two level
algebra of socle degree
$j$. For all pairs $(u,i=j-u)$ satisfying $0< u<j$
we have
\begin{equation}\label{eqmain}
H(A)_u-d_i(F,G) \le H(A_{\lambda_{gen}})_{i}.
\end{equation}
\end{uthm}\noindent
In Theorem \ref{mainthm2} we give a lower bound for 
$H(A_{\lambda_{gen}})$ that depends only on $H(A)$.
\par\smallskip\noindent
\subsection{Questions and examples: pencils of forms}\label{subsquest}
We offer some questions about
pencils of forms and the Hilbert functions they determine, and state their status.
This provides some further context for our work, and as well we pose
open problems.
\begin{question}\label{q1} What are
natural invariants for pencils of forms?
\begin{enumerate}[i.]
\item\label{q1i}What sequences $H$ occur as Hilbert functions
$H(A)$?\par
Status: Open for $r\ge 3$, even for $t=2$, but see \cite{GHMS,Za1}.
\item\label{q1ii} Is there a sequence $H=(1,3,\ldots ,2,0)$, such that
LevAlg(H) has two irreducible components?\par
Status:  Open. The answer to the analogous question is "no" for embedding
dimension $r=2$, and "yes" for $r\ge 4$.
\item\label{q1iii}Can we use our knowledge of the Hilbert functions and
parameter spaces for Artinian Gorenstein algebras, to study type two
level algebras $A$ of embedding dimensions three and four?
\par
Status: This has been the main approach to
classifying type two level algebras. See \cite{GHMS} and \cite{Za1,Za2}, as well
as Lemma \ref{lem1.1}, Examples
\ref{ex2.7} and \ref{ex2.8} below. 
\item\label{q1iv}Given a type two Artinian algebra $A$, consider the
pencil of Gorenstein quotients $A_\lambda = R/\Ann (F+\lambda G)$
having the same socle degree as $A$. What can be said about
the Hilbert functions $H(A_\lambda)$?\par
Status:  We begin a study here.
See also \cite{I4,Za1,Za2}.
\end{enumerate}
\end{question}
 The Question \ref{q1} about natural invariants of $A$ connects also
with classical invariant theory, but we do not pursue this here: see
\cite{DK} and \cite{RS} for analogous connections in the Gorenstein
case. The following example illustrates Question \ref{q1}\eqref{q1iv}, and as well
the main result.
\begin{example} Let $r=2$. $F=X^4,G=XY^3$, then
$H_F=(1,1,1,1,1)$, $H_G=(1,2,2,2,1)$, the ideal
 $I=\Ann(F,G)=\Ann F\cap \Ann G=
(x^2y,y^4,x^5)$. The type two level algebra $A=R/I$ has Hilbert
function
$ H(A)=(1,2,3,3,2)$. The dualizing module
$\hat{A}=R\circ \langle F,G\rangle \subset \mathcal D$ satisfies
\begin{equation*}
\hat{A}=\langle 1;X,Y;\ X^2,YX,Y^2;\ X^3,Y^3,
XY^2;\ X^4,XY^3\rangle.
\end{equation*}
The Gorenstein quotients $A_\lambda$ satisfy
\begin{equation*}
H(A_\lambda )=(1,2,3,2,1) \text { for } \lambda \not= \infty,0.
\end{equation*}
This is
the maximum possible (so compressed) Hilbert function for a
Gorenstein Artinian quotient of $R$ having socle degree 4.
\end{example}
The following specific question arose from a discussion with A. Geramita
about the Hilbert functions possible for type two level algebras. It
was the starting point of our work here.
\begin{question}\label{QTony} Let $F,G$ be two degree-$j$ homogeneous
polynomials, elements of $ \mathcal{D}=K[X_1,\ldots ,X_r]$, such that
\begin{enumerate}[i.]
\item $F,G$ together have at least $2r-2$ linearly independent first
partial derivatives, and
\item $ F,G$ together involve all $ r$ variables: this is equivalent to
the
$(j-1)$-order partials of
$F,G$ spanning $\langle X_1,\ldots ,X_r\rangle$.
\end{enumerate}
Does some linear combination $F_\lambda= F+\lambda G$ have $r$ linearly
independent first partial derivatives? 
\end{question}
We answer this question positively in Corollary
\ref{deriv}. Here are two examples to illustrate.
\begin{example} Let  $r=3, j=4,  F=X^4+Y^4, G=(X+Y)^4+Z^4$.
Then the pencil $V=\langle F,G\rangle$ involves all three variables, and
these forms together have four linearly independent first partials
$X^3,Y^3,(X+Y)^3,Z^3$. For \emph{all} $\lambda \not=0$ the form
$F_\lambda = F+\lambda G$ has
three linearly  independent first partials.
\end{example}
\begin{example} Let  $r=3, j=4, F=XZ^3, G=YZ^3$, Then $V=\langle
F,G\rangle$ involves all three variables and these forms together have only
$3=2r-3$ linearly independent first partials. Each
$F_\lambda$ has only 2 linearly independent first partials.
Thus, the hypothesis in Question \ref{QTony} that $V$ has at least $2r-2$
linearly independent first partial derivatives is necessary for the desired
conclusion.
\end{example}
Our work here is focussed primarly on the following question.
\begin{question} Given Hilbert functions $H_F=H(R/\Ann F)$, and
$H_G=H(R/\Ann G)$ for two degree-$j$ forms $F,G\in \mathcal{D}$, or given 
$H(A),\  A=R/\Ann(F,G)$,  determine the possible  Hilbert functions $H(A_\lambda)$
for $A_\lambda =R/\Ann F_\lambda \, , F_\lambda =F+\lambda
G$? Are there {\it{numerical}} restrictions on
the generic value $H(A_{{\lambda}_{gen}})$?
\end{question}
It is easy to give a partial answer. Evidently, by \eqref{ebasicexact}, we
can't have two values of
$\lambda$ with
\begin{equation}
H(A_\lambda)_i< H(A)_i/2.
\end{equation}
We may conclude that small $H(A_\lambda)$ are rare, given $H(A)$! \par
In our main results
we show that if $H_F$ and $H_G$ are small in comparison with
$H(A)$, then $H(A_{\lambda_{gen}})$ is large
(Theorem \ref{mainthm}). We then show a lower bound for $H(A_{\lambda_{gen}})$
in terms of $H(A)$ (Theorem \ref{mainthm2}). Several examples illustrate the
results (see especially Examples \ref{ex2.7},\ref{ex2.8}). Example~\ref{exnomin}
gives a pencil $A_\lambda$ of Gorenstein Artinian quotients not having a minimum
Hilbert function; and Example~\ref{excompress} gives a compressed type two
Artinian level algebra $A$ such that $A_{\lambda_{gen}}$ is not compressed
Gorenstein.

\section{Hilbert Functions for pencils of forms}\label{sec2}
In this section we state and prove our main results.
 We first give an exact sequence relating $A_F,A_G$, and $A$.
We define $R$-module homomorphisms
\begin{align*}
\iota &: A \to R/\Ann F
\oplus R/\Ann G \quad\\
&\qquad \qquad \qquad \qquad \iota (f)=(f\mod
\Ann F ,\, -f\mod \Ann G )\\
\pi &:  R/\Ann F
\oplus R/\Ann G \to  R/(\Ann F +\Ann G):\quad\\
&\qquad \qquad \qquad \qquad\pi(a,b)= (a+b) \mod (\Ann F +\Ann G).
\end{align*}
\begin{lemma}\label{exactseq} Let $F,G\in {\mathcal D}_j$ determine a
type two level Artinian quotient $A=R/I, I=\Ann (F,G)$ of $R$. There is an
exact sequence of $R$-modules
\begin{equation}\label{ebasicexact}
 0\to A\xrightarrow{\iota }  R/\Ann F
\oplus R/\Ann G\xrightarrow{\pi}
 R/(\Ann F +\Ann G) \to 0,
\end{equation}
whose dual exact sequence is
\begin{equation}\label{ebasicexactd}
0\to \langle R\circ F\rangle \cap \langle R\circ G \rangle
\xrightarrow{\pi^\ast} \langle R\circ F\rangle \oplus \langle R\circ G
\rangle
\xrightarrow{\iota^\ast} R\circ
\langle F,G\rangle \to 0.
\end{equation}
\end{lemma}
\begin{proof} 
Since the duality between $R_j$ and ${\mathcal D}_j$ is
exact, we have 
\begin{equation}\label{edual1}
(\Ann F+\Ann G)^\perp= \langle R\circ F\rangle \cap
\langle R\circ G \rangle .
\end{equation}
Thus the two sequences are dual. Evidently $\iota$ is an inclusion and
$\pi$ is a surjection. The kernel of $\iota^\ast$ consists of pairs
$(h_1\circ F,h_2\circ G)$ such that $h_1\circ F-h_2\circ G=0$; this is
evidently the image of $\pi^\ast$, so the sequences are exact.
\end{proof}\par  
We let $J=\Ann F+\Ann G$: it depends of course upon the
choice of the pair
$(F,G)\in \langle F,G\rangle$. We
denote by $H(R/J)=(1,d_1,\ldots ,d_j)$ the Hilbert function
$H(R/J)$ where $d_i=d_i(F,G)$. Thus we have
from
\eqref{edual1}, that, letting $i=j-u$,
the integer $d_i$ measures the overlap in degree $i$ between the inverse systems
determined by $F$ and $G$:
\begin{align}\label{eqdi}
d_i&=\dim_K \langle R_u\circ F\rangle \cap
\langle R_u\circ G \rangle \\
&=\dim_K R_u\circ F+\dim_K R_u\circ G-H(A)_i\label{eqdi2}\\
&=(H_F)_i+(H_G)_i-H(A)_i.\label{eqdi3}
\end{align}
The equalities \eqref{eqdi2}, \eqref{eqdi3} are immediate from \eqref{eqdualb} and
\eqref{ebasicexactd}.
 We set, again letting $i=j-u$,
\begin{equation}\label{eqti}
t_i=\dim_K(((\Ann F)_u\circ G)\cap 
((\Ann G)_u\circ F)),
\end{equation}
where $t_i=t_i(F,G)$ depends on the pair $(F,G)$.
 Recall from Remark
\ref{remgen} that ``generic
$\lambda$'' refers to a suitable open dense set of $\lambda\in \mathbb P^1$, that
is, to all but a finite number of values of $\lambda$. We
denote by
$d_i,t_i$ the integers $d_i(F,G)$ and
$t_i(F,G)$ defined just above.
The 
following main result shows that a small overlap between $R_u\circ F$ and
$R_u\circ G$ implies a large value for
$H(A_{\lambda _{gen}})_i$.
\begin{theorem}\label{mainthm} Let $A=R/\Ann (F,G)$ be a type two level
algebra of socle degree
$j$. If $\lambda$ is generic, then for all pairs $(u,i=j-u)$ satisfying $0< u<j$
we have
\begin{equation}\label{eqmain}
H(A)_u-d_i \le H(A_\lambda)_{i}\le
H(A)_u-t_i.
\end{equation}
The upper bound on $H(A_\lambda)$ holds for all
$\lambda \ne 0, \infty$.
\end{theorem}
\begin{proof}
 Fix for now, and through the proof of the Claim below, an integer $u$ satisfying
$0<u<j$. By $``\dim V"$ below we mean $\dim_K V$. Let 
$C_u\subset R_u$ be a vector subspace complement to $(\Ann F)_u$, so
$C_u\oplus (\Ann F)_u=R_u$.
 Let $d=d_i$ and
let $e=\dim ((\Ann(F,G))_u)$: so $\dim A_u=\dim R_u-e$, and
let
$B\subset (\Ann F)_u$ be the vector subspace satisfying
\begin{equation}
B=\{ h\in (\Ann F)_u \mid h\circ G\in R_u\circ F\}.
\end{equation}
The homomorphism $h\to h\circ g,h\in B, g\in G$ induces a short exact sequence
\begin{equation*}
0\to (\Ann(F,G))_u\to B\to \langle R_u\circ F\rangle\cap \langle R_u\circ
G\rangle ,
\end{equation*}
implying
\begin{equation}
\dim B\le d+e.
\end{equation}
Since $(\Ann F)_u\circ (F+\lambda G)=(\Ann F)_u\circ G$, we have
\begin{equation}\label{eqkey}
R_u\circ (F+\lambda G)=C_u\circ
(F+\lambda G)+(\Ann F)_u\circ G.
\end{equation}
\begin{ucla} For generic $\lambda$ 
\begin{align}
\dim R_u\circ (F+\lambda G)&\ge \dim (C_u\circ F+ (\Ann F)_u\circ
G)\label{eqkey2}\\ &=\dim C_u+\dim (\Ann F)_u-\dim B\\
&=\dim R_u-\dim B\notag\\
&\ge\dim R_u-(d+e)\notag\\
&= \dim A_u-d,
\end{align} 
so $\dim R_u\circ (F+\lambda G)\ge \dim A_u-d$.
\end{ucla}
\begin{proof}[Proof of Claim]
The key step is \eqref{eqkey2}, which
results from \eqref{eqkey} and deformation. For the space $\langle
C_u\circ (F+\lambda G)+ (\Ann F)_u\circ G\rangle$ in \eqref{eqkey} is a
deformation of the space
$\langle C_u\circ F+ (\Ann F)_u\circ G\rangle$ on the right of
\eqref{eqkey2}, and dimension is a semicontinuous invariant. The other
steps are straightforward.
\end{proof}\par\noindent
The Claim shows the left hand inequality in \eqref{eqmain} for a
specific $u$.  Since $\mathbb
P^1_K$ is irreducible, the intersection of the dense open subsets of
$\mathbb P^1_K$ over which the left side of \eqref{eqmain} is satisfied
for each
$u, 0< u <j,$ is itself a dense open subset, completing the proof
that the
left side of \eqref{eqmain} holds simultaneously for generic $\lambda$ and
all such
$u$. \par
Suppose $\lambda \ne 0, \infty$. Let $C'_u$ be a complement in $R_u$ to
$J_u=(\Ann F)_u+(\Ann G)_u$. Then we have
\begin{align*}
R_u\circ (F+\lambda G)&=C'_u\circ (F+\lambda G)+(\Ann F)_u\circ G+(\Ann
G)_u\circ F\text { implying }\\
\dim R_u\circ (F+\lambda G)&\le \dim_k R_u-\dim(\Ann F)_u\cap
(\Ann G)_u)-t_i  =\dim A_u -t_i.
\end{align*}
This completes the proof of Theorem \ref{mainthm}. 
\end{proof}  
\begin{example}[No minimum $H(A_\lambda)$]\label{exnomin}
 Let $r=3$, $G=X^8+Y^4Z^4$, and $
F=L_1^8+\cdots +L_5^8$, where the
$L_i=a_{i1}X+a_{i2}Y+a_{i3}Z$ are general enough linear forms, elements of
$\mathcal{D}_1$. Here ``general enough'' means that their coefficients $\{
a_{ij}\in K\}$ lie in the open dense subset of the affine space $\mathbb
A^{15}$ such that the powers $L_1^{j-u},\ldots
,L_5^{j-u}$ are linearly independent in
$\mathcal{D}_{j-u}$ and maximally disjoint from $R_u\circ G$, $2\le u\le 6$ (see
\cite{I2}). Then we have for
$3\le u\le 6$
\begin{equation}
R_u\circ F=\langle L_1^{j-u},\ldots ,L_5^{j-u}\rangle,
\end{equation}
satisfying $\dim_k R_u\circ F=5$. This determines $H_F$, and we have
\begin{align*} 
H_F&= (1,3,5,5,5,5, 5,3,1)\\
H_G&= (1,3,4,5,6,5,4,3,1)\\
H(A)&=(1,3,6,10,11,10,9,6,2)=H_F+_h H_G,
\end{align*}
where by $H_F+_hH_G$ we mean the sequence satisfying
\begin{equation}\label{eqhadd}
(H_F+_hH_G)_i=\min\{\dim R_i, (H_F)_i+(H_G)_i\}.
\end{equation}
Theorem 1 implies that 
$H(A_{\lambda_{gen}})=(1,3,6,10,11,10,6,3,1)$. It is easy to check that
there are no other values of $\lambda$ other than $0, \infty$
(corresponding to $F,G$) such that
$H(A_\lambda)$ is smaller than $H(A_{\lambda_{gen}})$, and hence  no
minimum sequence $H(A_\lambda)$, since $H_F$ and $H_G$ are incomparable. 
\end{example}
Our second main result gives a lower bound for $H(A_{\lambda_{gen}})$ solely in
terms of $H(A)$.
\begin{theorem}\label{mainthm2}  Let $A$ be a type two level Artinian algebra of
socle degree $j$, and let
$u,i$ satisfy
$0<u\le i=j-u$. Assume that
$ H(A)_i\ge 2H(A)_u-2-3\delta _u$, where $\delta_u\ge 0$ and $\delta_u$ is an
integer. Then 
\begin{equation}\label{eq2.15}
H(A_{{\lambda}_{gen}})_i  \ge H(A)_u-\delta_u.
\end{equation}
\end{theorem}
\begin{proof} Assume the hypotheses of the Theorem, and suppose by way of
contradiction, that for generic $\lambda$ there is an integer $a\ge 0$ satisfying
\begin{equation}\label{eqast}
H(A_\lambda)_i=H(A)_u-\delta_u-1-a.
\end{equation}
 Take two generic forms $F',G'$ in the pencil.
Then the overlap between
$R_u\circ F'$ and $ R_u\circ G'$ (see \eqref{eqdi},\eqref{eqdi2})
satisfies
\begin{align}
d_i&=2(H(A)_u-\delta_u-1-a)-H(A)_i\notag\\
&\le 2H(A)_u-2\delta_u-2-2a-(2H(A)_u-2-3\delta _u) \notag\\
&\le \delta_u-2a.\label{eqastast}
\end{align}
By Theorem \ref{mainthm}, for generic $\lambda$ the AG quotient
$A'_\lambda=R/\Ann(F'+\lambda G')$ satisfies 
\begin{equation}
H(A'_\lambda)_i\ge  H(A)_u-(\delta_u-2a),
\end{equation}
 a
contradiction with equation \eqref{eqast}.
It follows that the assumed equation \eqref{eqast} is false, hence
\begin{equation}
H(A_\lambda)_i\ge H(A)_u-\delta_u
\end{equation}
 which is Theorem \ref{mainthm2}.
\end{proof}\par
\begin{unote}
 We assumed in the statement and proof of Theorem \ref{mainthm2} that
$\delta_u$ is an integer. Alternatively we could define $\delta
'_u=(2H(A)_u-2-H(A)_i)/3$ and conclude that $H(A_{{\lambda}_{gen}})_i \ge
H(A)_u-\lceil \delta'_u\rceil$ when $\delta'_u\ge 0$, and
$H(A_{{\lambda}_{gen}})_i\ge H(A)_u$ otherwise.
\end{unote}
 In the following corollary we give
a positive answer to Question
\ref{QTony}.
\begin{corollary}\label{deriv} Let $F,G$ together have at least $2r-2$ linearly
independent first partial derivatives, and suppose that
$F,G$ involve all $r$ variables. Then $\dim_K R_1\circ (F+\lambda G)=r$ for generic
$\lambda$.
\end{corollary}
\begin{proof} Take $i=j-1,u=1, \delta_1=0$ in Theorem \ref{mainthm2}.
From the assumptions we have
\begin{equation*}
H(A)_{j-1}\ge 2r-2=2 H(A)_1-2,
\end{equation*} 
 which implies by Theorem
\ref{mainthm2} that for generic $\lambda$ the dimension
 $H(A_\lambda)_{j-1}=r.$
\end{proof}\par
In order to apply Theorem \ref{mainthm2} most effectively, we use the following
result from \cite{GHMS}. For a sequence
$H=(1,\ldots ,H_j),\ H_j>0$ we denote by $H\ \hat{}\ $ the reverse sequence ${H\
\hat{}}_i=H_{j-i}, 1\le i\le j$. Recall that an
$O$-sequence is one that is the Hilbert function of some Artinian algebra
\cite{Mac2,BH}.
\begin{lemma}\label{lem1.1} (A. Geramita et al \cite{GHMS}) Let $A$ be a type $t$
level algebra with dualizing module $\hat{A}$. Let $A=R/\Ann W, W=W_A\subset
{\mathcal D}_j, \dim_K W=t\ge 2$ and let $V\subset W$ be a vector subspace of
codimension one. Then there is an exact sequence of
$R$-modules relating the type $t-1$ level algebra $B_V=B/\Ann V$ to $A$,
\begin{equation}
   0\to C\to A \to B_V \to 0
\end{equation}
whose dual exact sequence of $R$ submodules of $\mathcal D$ is
\begin{equation}\label{exseq}
 0\to \hat{B_V} \to \hat{A}\to \hat{C}\to 0.
\end{equation}
Here $\hat{C}$ is a simple $R$-submodule (single generator).
 We have for their Hilbert functions
\begin{equation*}
H(A)=H(B_V)+H(C), 
\end{equation*}
and the reverse sequence $H(C)\hat{}=(1,\ldots )$ is an
$O$-sequence.
\end{lemma}
\begin{proof} Let $F\in V$ span a complement to $V$ in $W_A$. Let the homomorphism
\begin{equation*}
 \tau _{(F,W,V)}:  R\to  R\circ F/\langle R\circ W\cap R\circ F\rangle\to 0
\end{equation*}
have kernel $S$. The module $\hat{C}$ is isomorphic to $R/S$. This shows
\eqref{exseq} and that $\hat{C}$ is simple.
\end{proof}
\begin{example}\label{ex2.7}  Let $r=3$ and suppose $A$ is a type two level
algebra satisfying
$H(A)=(1,3,\ldots ,4,2)$. Then the pencil $W_A=\langle F,G\rangle\subset \mathcal
D_j$ defining $A$ may be chosen so that $H_F=(1,3,\ldots ,3,1)$ and 
\begin{equation}\label{eqspecialsum}
H(A)\le
(1,3,6,10,\ldots ,6,3,1)+_h(1,1,\ldots ,1),
\end{equation}
(see \eqref{eqhadd} for the sum $+_h$ used above).
In particular $(1,3,\ldots ,8,4,2)$ and $(1,3,\ldots ,12,7,4,2)$ are not
sequences possible for the Hilbert function of a level algebra quotient of
$R=K[x,y,z]$. Here in equation \eqref{eqspecialsum}, the sequence $(1,3,6,\ldots
,3,1)$ is the compressed Gorenstein sequence of socle degree
$j$ (see below).\par 
Here the Corollary \ref{deriv} implies that we may choose $G\in W_A$
such that $H_G=(1,3,\ldots ,3,1)$; it follows that
$H(\hat{C})=(1,1,\ldots )$, so to be an $O$-sequence, $H(\hat{C})\le
(1,1,\ldots ,1)$; this and \eqref{exseq} show
\eqref{eqspecialsum}. F. Zanello has extended this kind
of result, and in \cite{Za1} shows sharp upper bounds for the Hilbert function $H(A)$
for type two level algebra quotients of $R$ in $r$-variables given $H(A)_{j-1}$.
\end{example}
\begin{example}\label{ex2.8} Let $r=3$ and let $H=(1,3,6,8,6,4,2)$. Consider
first
$W_1=\langle F,G\rangle$ where $H_F=(1,3,5,7,5,3,1)$, and $G=L^6, L$ a general
enough linear form (element of ${\mathcal D}_1$), so $H_G=(1,1,\ldots ,1)$.
Then
$A_1=R/\Ann W_1$ is easily shown to have Hilbert function $H$, as, choosing $F$
first, and $G$ second,
$R_u\circ G=\langle L^{6-u}\rangle$ and is linearly disjoint in general from
$R_u\circ G$, by the spanning property of the rational normal curve: powers of
linear forms span $D_i,i=j-u$ (see
\cite{I2}).
\par Next, let $W_2=\langle F',G'\rangle$ where $F'$ is a general element of
$K[X,Y]_6$, and $G'$ is a general enough element of $K[X+Y,Z]_6$. Then
$H_F=H_G=(1,2,3,4,3,2,1)$ and $A_2=R/\Ann W_2$  also has Hilbert function $H$.
In either case, Theorem \ref{mainthm} implies that 
\begin{equation*}
H(R/\Ann F_{\lambda_{gen}})=
H(R/\Ann F'_{\lambda_{gen}})=(1,3,6,8,6,3,1).
\end{equation*}
 In this example, it is the 
non-generic elements of each pencil --- the ``unexpected properties'' --- that serve to
distinguish the pencil.
\end{example}
 A \emph{compressed} level algebra of given type $t$, socle degree $j$, and
embedding dimension $r$ is one having the maximum possible Hilbert function
given those integers (see \cite{I1,FL,Bj1,MMN}).\par
The following example responds negatively in embedding dimension three to a
question  of D.~Eisenbud, B.
Ulrich, and C.~Huneke. They asked if a generic socle degree-$j$
AG quotient $A_{\lambda_{gen}}$ of a compressed type two level
 algebra $A$,  must also be compressed. This is true for $r=2$.
\begin{example}\label{excompress} Let
$r=3,j=4$, take $a,b\in K, (a,b)\not=(0,0)$, and set 
\begin{equation*}
 F=X^3Y+X^2Z^2+aXZ^3+bYZ^3, \quad
G=X^3Z+X^2Y^2+X^2YZ+3aXY^2Z+bY^3Z.
\end{equation*}
Here $H(A)=(1,3,6,6,2)$, and is compressed,  but $(y^2-\lambda z^2)\circ
(F+\lambda G)=0$ and we have that for generic $\lambda$, $H(A_\lambda ) =
(1,3,5,3,1)$,
 which is not the compressed sequence $(1,3,6,3,1)$. Note that in
applying Theorem
\ref{mainthm2} with $i=2,j=4$ here we would have
$6=H(A)_2=2H(A)_{4-2}-2-3(4/3)$, so taking 
$\delta =2$ in equation \eqref{eq2.15}, we would conclude only that 
$H(A_{\lambda_{gen}})_2\ge 6-2=4$.
\end{example}
\begin{remark}[Characteristic $p$]\label{charp} Assume that $K$ is an
infinite field of finite characteristic $p$. For $p>j$, the socle degree,
there is no difference in statements. For $p\le j$ one must use the
divided power ring in place of $\mathcal D$, and the action of $R$ on
$\mathcal D$ is the contraction action. With that substitution, the lemmas and
theorems here extend to characteristic $p$. However, in examples, one must
substitute divided powers for regular powers --- see 
\cite[Appendix A]{IK} for further discussion. 

\end{remark}

\begin{ack}
This article began after a conversation with Tony Geramita, in Fall 2002
during a visit to the MSRI Commutative Algebra Year. We noted that
in tables he and colleagues had calculated for Hilbert functions of type
two algebras (see
\cite{GHMS}), that Corollary
\ref{deriv} was satisfied. I believed there would be a general
result, which turned out to be Theorem~\ref{mainthm}.
I thank Tony for this discussion, and for helpful comments. I am appreciative to the
organizers of the Siena conference ``Projective Varieties with Unexpected Properties''
for hosting a lively and informative meeting, and the impetus to write a concise
presentation of the results; and I thank the referee.
\end{ack}

\bibliographystyle{amsalpha} 

\end{document}